\documentclass[a4paper,11pt]{article}

\usepackage{fullpage}
\usepackage[dvips]{epsfig}
\usepackage{psfrag}
\usepackage{amsfonts}
\usepackage{latexsym}
\usepackage{amsmath}
\usepackage{amssymb}
\usepackage{mathrsfs}
\usepackage{color}
\usepackage{theorem}
\usepackage{soul}

{\theoremstyle{plain}
\newtheorem{theo}{Theorem}
\newtheorem{prop}[theo]{Proposition}
\newtheorem{cor}[theo]{Corollary}
\newtheorem{lem}[theo]{Lemma}

{\theorembodyfont {\rmfamily}
\newtheorem{exa}[theo]{Example}
\newtheorem{nota}[theo]{Notation}
\newtheorem{defn}[theo]{Definition}
\newtheorem{rem}[theo]{Remark}

}}

\begin{document}

\title{On a generalization of Christoffel words: epichristoffel words\thanks{with the support of {\it Fonds de recherche sur la nature et les technologies} (Qu\'ebec, Canada)}}
\author{Genevi\`eve Paquin\\
Laboratoire de math\'ematiques,  Universit\'e de Savoie, CNRS UMR 5127\\
73376 Le Bourget du Lac, France\\
{\tt Genevieve.Paquin@univ-savoie.fr}}
%

\date{\today}

\maketitle
\sloppy

\newcommand{\QED}{\rule{1ex}{1ex} \par\medskip}

\newcommand{\nl}{\par\medskip\noindent}
\newcommand{\Proof}{{\nl\it Proof.\ }}
\newcommand{\Ind}{\hbox{\rm Ind}}
\newcommand{\Predec}{\hbox{\rm pred}}
\newcommand{\Succes}{\hbox{\rm succ}}
\newcommand{\Label}{\hbox{\rm label}}
\newcommand{\card}{\hbox{\rm Card}}

\newcommand{\Pref}{\hbox{\rm Pref}}
\newcommand{\Suff}{\hbox{\rm Suff}}
\newcommand{\Pal}{\hbox{\rm Pal}}
\newcommand{\Squares}{\hbox{\rm Squares}}
\newcommand{\Overlaps}{\hbox{\rm Overlaps}}
\newcommand{\Lpc}{\hbox{\rm Lpc}}
\newcommand{\Card}{\hbox{\rm Card}}
\newcommand{\Last}{\hbox{\rm Last}}
\newcommand{\First}{\hbox{\rm First}}
\newcommand{\fleche}{\longrightarrow}
\newcommand{\donne}{\longmapsto}
\newcommand{\vphi}{\varphi}
\newcommand{\N}{\mathbb N}
\newcommand{\Z}{\mathbb Z}
\newcommand{\bprop}{\begin{prop}}
\newcommand{\eprop}{\end{prop}}
\newcommand{\bcor}{\begin{cor}}
\newcommand{\ecor}{\end{cor}}
\newcommand{\blem}{\begin{lem}}
\newcommand{\elem}{\end{lem}}

\newcommand{\A}{\mathcal{A}}
\newcommand{\R}{\mathbb{R}}
\newcommand{\Q}{\mathbb{Q}}
\newcommand{\Ult}{\textnormal{Ult}}
\newcommand{\al}{\textnormal{Alph}}

\begin{abstract}  Sturmian sequences are well-known as the ones having minimal complexity over a 2-letter alphabet. They are also the balanced sequences over a 2-letter alphabet and the sequences describing discrete lines. They are famous and have been extensively studied since the 18th century. One of the {extensions} of these sequences over a $k$-letter alphabet, with $k\geq 3$, are the episturmian sequences, which generalizes a construction of Sturmian sequences using the palindromic closure operation. There exists a finite version of the Sturmian sequences called the Christoffel words. They are known since the works of Christoffel  and have interested many mathematicians. In this paper, we introduce a generalization of Christoffel words  for an alphabet with 3 letters or more, using the episturmian morphisms. We call them the {\it epichristoffel words}. We define this new class of finite words and show how some of the properties of the Christoffel words can be generalized naturally or not for this class.  
\end{abstract}


\section{Introduction}

As far as we know, Sturmian sequences first appeared  in the literature at the 18th century in the precursory works of the astronomer Bernoulli \cite{jb1772}. They later appeared in the 19th century in Christoffel \cite{ebc1875} and Markov \cite{am1882} works. The first deep study of these sequences is given in \cite{mh1938,mh1940} where the name {\it Sturmian sequence} appears for the first time. At the end of the 20th century and more recently, many mathematicians have been interested in those sequences, for instance \cite{ch1973,emc1974,kbs1976,tcb1993,bpr1994,gz1995,adl1997,jb2002}.  Recent books also show this interest \cite{ml2002,pf2002,as2003,blrs2008} as well as a recent survey \cite{jb2007}. In this wide literature, we find different characterizations of the Sturmian sequences. In particular, they are the sequences over a 2-letter alphabet having the minimal complexity, they also are the balanced sequences over a 2-letter alphabet and they code discrete lines. These different characterizations show how the Sturmian sequences occur in different fields as number theory \cite{rm19852,rjs1991,rt2000,rt2001,bv2003,rjs2004,go2005}, discrete geometry, crystallography \cite{bt1986} and symbolic dynamics \cite{mh1938,mh1940,gah1944,mq1987}. 

Since the end of the 20th century, numerous generalizations of Sturmian sequences have been introduced for an alphabet with more than 2 letters. Among them, one natural generalization is called the {\it episturmian sequences} and is using the palindromic closure property of Sturmian sequences \cite{adl19972}. The first construction of episturmian sequences is due to \cite{djp2001}.  Previously the first introduction and study of an episturmian sequence was that of the Tribonacci word \cite{gr1982} and an important class of episturmian sequences, now called the Arnoux-Rauzy sequences, had been considered in \cite{ar1991,rz2000}. More recently the whole class was extensively studied, for instance in \cite{jv2000,rz2000,djp2001,jp2002,jp2004,jj2005,ag2005,pv2007,gr20072,ag2006,bdldlz2008,gjp2007,glr2008}.  For  surveys about episturmian sequences, see for instance \cite{jb2007,gj2009}.

The finite version of Sturmian sequences, called  {\it Christoffel words}, has been also well studied  \cite{ebc1875,ml2002,br2006,bdlr2007,kr2007}. It is known that any finite standard Sturmian word, that is the words obtained by standard Sturmian morphisms to a letter, is conjugate to a Christoffel word. A Christoffel word is then the smallest word, with respect to the lexicographic order, in the conjugacy class of a finite standard Sturmian word.  Finite factors of the episturmian sequences appeared for instance in \cite{gjp2007}. The class of standard episturmian words is naturally defined as the set of finite words obtained by standard episturmian morphisms to letter, but no generalization of the Christoffel words have been introduced yet. In this paper, we introduce such a generalization that we naturally call the {\it epichristoffel words}. Note that it naturally appears that for each standard episturmian word, there exists a conjugate which is an epichristoffel word, and conversely. 

The paper is organized as follows.

We first recall some basic definitions of combinatorics on words and we establish the notation used in this paper. We recall the definitions and some properties of the Sturmian sequences,  the Christoffel words and the episturmian sequences. Then we introduce our new class of finite words: the {\it epichristoffel ones}.  We prove how some of the properties of the Christoffel words can be generalized for an alphabet with more than $2$ letters.   We then describe an algorithm which determines if a given $k$-tuple describes the occurrence numbers of letters in an epichristoffel word or not. If so, we show how to construct it. Finally, we prove the next theorem, which is a generalization  of a result for Christoffel words \cite{dldl2006}, that characterizes epichristoffel conjugates. \\

\noindent {\bf Theorem }   {\it Let $w$ be a finite primitive word different from a letter. Then the conjugates of $w$ are all factors of the same episturmian sequence if and only if  $w$ is conjugate to an epichristoffel word.
}


\section{Definitions and notation}

Throughout this paper, $\A$ denotes a finite alphabet containing  $k$ letters $a_0, a_1, \ldots, a_{k-1}$. A {\it finite word} is an element of the free monoid $\A ^*$. If $w=w[0]w[1]\cdots w[n-1]$, with $w[i] \in \A$,  then $w$ is said to be a finite word of {\it length} $n$ and we write $|w|=n$. By convention, the {\it empty word} is denoted $\varepsilon$ and its length is 0. We define $\A ^\omega$ the set of right infinite words, also called {\it sequences}, over the alphabet $\A$ and then, $\A ^\infty=\A^* \cup \A^\omega$ is the set of finite and infinite words. 

The number of occurrences of the letter $a_i$ in $w$ is denoted $|w|_{a_i}$.  The {\it reversal} of the word $w=w[0]w[1]\cdots w[n-1]$ is $\widetilde{w}=w[n-1]w[n-2] \cdots w[0]$ and  if $\widetilde{w}=w$,  then $w$ is said to be a {\it palindrome}. A finite word $f$ is a {\it factor} of $w \in \A^\infty$ if $w=pfs$ for some $p \in \A^*, s \in \A^\infty$. If $p=\varepsilon$ (resp. $s=\varepsilon$), $f$ is called a {\it prefix} (resp. a {\it suffix}) of $w$.  Let $u\in \A^*$ and $n\in \N$. We denote by $u^n$ the word $u$ repeated $n$ times and we called it a {\it $n$-th power word}. A factor $\alpha^k$ of the word $w$, with $\alpha \in \A$ and $k \in \N$ locally maximum, is called a {\it block of $\alpha$ of length $k$ in $w$}. Let $u, v$ be two palindromes, then $u$ is a {\it central factor} of $v$ if $v=wu\tilde w$ for some $w \in \A^*$. The {\it right palindromic  closure} of $w \in \A^*$ is the shortest palindrome $u=w^{(+)}$ having $w$ as prefix. 

The {\it set of factors} of $w \in \A^\omega$ is denoted $F(w)$ and $F_n(w)=F(w)\cap \A ^n$ is the set of all factors of $w$ of length $n \in \N$.  The complexity function is given by $P(n)=|F_n(w)|$ and is the number of distinct factors of $w$ of length $n \in \N$. Two words $w$ and $w'$ are said {\it equivalent} if  they have the same set of factors: $F(w)=F(w')$.

The {\it conjugacy class} $[w]$ of $w \in \A ^n$ is the set of all words $w[i]w[i+1]\cdots w[n-1]w[0]\cdots w[i-1]$, for $0 \leq i \leq n-1$. If $w$ is not the power of a shorter word, then $w$ is said to be  {\it primitive} and has exactly $n$ conjugates. If $w$ is the smallest of its conjugacy class, relatively to some lexicographic order, then $w$ is called {\it a Lyndon word}. 

Let $w$ be an infinite word, then a factor $f$ of $w$ is {\it right} (resp. {\it left}) {\it special} in $w$ if there exist $a,b \in \A$, $a \neq b$, such that $fa, fb \in F(w)$ (resp. $af, bf \in F(w)$).  A word $w$ over $\A$ is {\it balanced} if for all factors $u$ and $v$ of $w$ having the same length, for all letters $a \in \A$, one has
$$\big ||u|_a-|v|_a\big |\leq 1.$$
 
 If $w=pus \in \A^\infty$, with $p,u \in \A^*$ and $s \in \A^\infty$, then $p^{-1}w$ denotes the word $us$. Similarly, $ws^{-1}$ denotes the word $pu$.
 
An integer $p\in \N$ is a {\it period} of the word $w=w[0]w[1]\cdots w[n-1]\in \A^*$ if $w[i]=w[i+p]$ for $0\leq i < n-p$.  When $p=0$, the period is {\it trivial}. If $p$ is the smallest non trivial period of $w$, then the {\it fractionnary root} of $w$ is defined as the  prefix  $z_w$ of $w$ of length $p$. An infinite word $w\in \A^\omega$ is {\it periodic} (resp. {\it ultimately periodic}) if it can be written as $w=u^\omega$ (resp. $w=vu^\omega$), for some $u,v \in \A^*$.  If $w$ is not ultimately periodic, then it is {\it aperiodic}. 
A {\it morphism} $f$ from $\A^*$ to $\A^*$ is a mapping from $\A^*$ to $\A^*$ such that for all words $u,v\in\A^*$, $f(uv)=f(u)f(v)$.  A morphism extends naturally on infinite words.

\section{Sturmian, Christoffel and episturmian words}

Before introducing our generalization of  Christoffel words, inspired by the definition of episturmian sequences, let us recall the definition of these well-known families and some of their properties.

\subsection{Sturmian words and morphisms}


One of the classical definitions of Sturmian sequences is the one given by Morse and Hedlund \cite{mh1940}: \\
	
\noindent {\bf Definition } Let $\rho$, called the intercept, and $\alpha$, called the slope, be two real numbers with $\alpha$ irrational such that $0\leq \alpha< 1$. For $n\geq 0$, let 
\begin{center}
$s[n] = \left \{ \begin{tabular}{l}
$a$ \textnormal{if} $\lfloor \alpha(n+1) +\rho  \rfloor = \lfloor \alpha n +\rho \rfloor$,\\
$b$ \textnormal {otherwise},
\end{tabular} \right . $ \\ 

$s'[n] = \left \{ \begin{tabular}{l}
$a$ \textnormal{if} $\lceil \alpha(n+1) +\rho  \rceil = \lceil \alpha n +\rho \rceil$,\\
$b$ \textnormal {otherwise}.
\end{tabular} \right . $\\
\end{center}
Then the sequences 
$$s_{\alpha,\rho}=s[0] s[1]s[2]\cdots \quad \quad \textnormal{ and } \quad \quad s'_{\alpha, \rho}=s'[0]s'[1]s'[2]\cdots$$ are 
{\it Sturmian} and conversely, a Sturmian sequence can be written as $s_{\alpha, \rho}$ or $s'_{\alpha,\rho}$ for $\alpha$ irrational and $\rho \in \R$. \\

Sturmian sequences have several characterizations. For more details about this class of words,  we refer the reader to the section in \cite{ml2002} devoted to Sturmian sequences.\\

\noindent {\bf Proposition }  \cite{ch1973}  A sequence $s$ is Sturmian if and only if  for all $n \in \N$, $P(n)=n+1$.\\

\newpage
\begin{theo}{\rm \cite{ml2002}} Let $s$ be a sequence. The following assertions are equivalent: \vspace{-0.2cm}
\begin{itemize}
\item [\rm i)] $s$ is Sturmian;
\item [\rm ii)] $s$ is balanced and  aperiodic.
\end{itemize}
\end{theo}

\begin{defn}\cite{ml2002} A {\it morphism} $f$ is  {\it Sturmian} if  $f(s)$ is Sturmian for all Sturmian sequences~$s$. 
\end{defn}


\subsection{Christoffel words}

In discrete geometry, Christoffel words are defined as the discretization of a line having a rational slope, as introduced in  \cite{bl1993}. In symbolic dynamics, they are defined by exchange of intervals  \cite{mh1940} as follows.\\

\noindent {\bf Definition }  Let $p$ and $q$ be positive relatively prime integers and $n=p+q$. Given an ordered $2$-letter alphabet $\{a < b\}$, the {\it Christoffel word w of slope $p/q$} over this alphabet is defined as $w=w[0]w[1]\cdots w[n-1]$, with 
$$w[i] = \left \{ \begin{tabular}{l}
$a$ \textnormal{if} $ip \mod n > (i-1)p \mod n,$ \\
$b$ \textnormal{if} $ip \mod n < (i-1)p \mod n,$
\end{tabular}
\right . $$ for $0\leq i\leq n-1$, where $k \mod n$ denotes the remainder of the Euclidean division of $k$ by $n$.\\

Notice that since $p$ and $q$ are relatively prime, a Christoffel word is always primitive. {Other important properties of Christoffel words will be recalled just before their generalizations in Section 4.}

\subsection{Episturmian sequences and morphisms}

One of the possible generalizations of  Sturmian sequences for an alphabet with $3$ letters or more is the set of episturmian sequences. Let us first recall the definition of standard episturmian sequences as introduced initially  by Droubay, Justin and Pirillo.

\begin{defn} \cite{djp2001} \label{def1epis} A sequence $s$ is  {\it standard episturmian} if it satisfies one of the following equivalent conditions.
\vspace{-0.2cm}
\begin{itemize}
\item [\rm i)] For every prefix $u$ of $s$, $u^{(+)}$ is also a prefix of  $s$.
\item [\rm ii)] Every leftmost occurrence of a palindrome in $s$ is a central factor of a  palindromic prefix  of $s$.
\item [\rm iii)] There exist a sequence $u_0=\varepsilon, u_1, u_2, \ldots$ of palindromes and a sequence  $\Delta(s)=x[0]x[1]\cdots$, with $x[i] \in \A$, such that $u_n$ defined by $u_{n+1}=(u_nx[n])^{(+)}$, with $n \geq 0$,  is a prefix of $s$.
\end{itemize}
\end{defn}

\begin{defn} \cite{djp2001} \label{episst} A sequence $t$ is {\it episturmian} if $F(t)=F(s)$ for a standard episturmian sequence $s$.   
 \end{defn}

An equivalent definition is that a sequence $s \in A^\omega$ is {\it episturmian} if its set of factors is closed under reversal and $s$ has at most one right (or equivalently left) special factor for each length.

\begin{nota} \textnormal{\cite{jj2005}} Let  $w=w[0]w[1]\cdots w[n-1]$, with $w[i] \in \A$, and $u_0=\varepsilon$, \ldots, $u_{n}=(u_{n-1}w[n-1])^{(+)}$, the palindromic prefixes of  $u_{n}$. Then \textnormal{Pal}$(w)$\index{$\textnormal{Pal}(w)$} denotes the word $u_{n}$. 
\end{nota}

In Definition \ref{def1epis}, $\Delta(s)$ is called the {\it  directive sequence} of the standard episturmian sequence~$s$. Since $\Delta(s)$ is the limit of its prefixes and $s$ is the limit of the $u_n$, it is natural to write $s=\Pal(\Delta(s))$.

Let us recall  from   \cite {jj2005} a useful property of the operator $\Pal$.

\blem \label{lemjj} {\rm \cite{jj2005}} Let $x \in \A$, $w \in \A^*$. If $|w|_x=0$, then $\Pal(wx)=\Pal(w)x\Pal(w)$. Otherwise, we write $w=w_1xw_2$ with $|w_2|_x=0$. The longest palindromic prefix of  $\Pal(w)$ which is followed by  $x$ in $\Pal(w)$ is $\Pal(w_1)$. Thus, $\Pal(wx)=\Pal(w)\Pal(w_1)^{-1}\Pal(w)$. 
\elem

\begin{defn} For $a, b \in \A$, we define the following  endomorphisms of $\A ^*$:
\vspace{-0.2cm}
\begin{itemize}
\item [\rm i)] $\psi_a(a)=\overline{\psi }_a(a)=a$;
\item [\rm ii)]  $\psi_a(x)=ax$, if $x \in \A \setminus \{a\}$;
\item [\rm iii)] $\overline{\psi }_a(x)=xa$, if $x \in \A \setminus \{a\}$;
\item [\rm iv)] $\theta_{ab}(a)=b$ , $\theta_{ab}(b)=a$, $\theta_{ab}(x)=x$, $x \in \A \setminus \{a,b\}$. 
\end{itemize}
\end{defn}

The endomorphisms $\psi$ and $\overline \psi$ can be naturally extended to a finite word $w=w[0]w[1]\cdots w[n-1]$. Then $\psi_w(a)=\psi_{w[0]}(\psi_{w[1]}(\cdots (\psi_{w[n-1]}(a))\cdots ))$ and $\overline \psi_w(a)=\overline \psi_{w[0]}(\overline \psi_{w[1]}(\cdots (\overline \psi_{w[n-1]}(a))\cdots ))$ , with $a \in \A$.

Similarly to the Sturmian morphisms, we can define the episturmian morphisms as follows.

\begin{defn}\cite{jp2002} The set $\mathscr{E}$ of {\it episturmian morphisms} is the monoid generated by the morphisms $\psi_a, \overline{\psi}_a, \theta_{ab}$ \index{$\psi_a$}\index{$\overline{\psi}_a$}\index{$ \theta_{ab}$} under composition. The set  $\mathscr{S}$ of {\it standard episturmian morphisms}  is the submonoid generated by the $\psi_a$ and  $\theta_{ab}$; the set of {\it pure episturmian morphisms} is the submonoid generated by the $\psi_a$ and $\overline \psi_a$.
\end{defn}

As the Sturmian morphism, the episturmian ones have the following characteristic property: a morphism $f$ is {\it episturmian} if $f(s)$ is episturmian for any episturmian sequence~$s$.


\section{Epichristoffel words}

In this section, we generalize  Christoffel words to a $k$-letter alphabet and we call this generalization {\it epichristoffel words}. 

Let us first recall some properties of Christoffel words that will be used to define their generalization.

\blem {\rm \cite{bdl1997}} \label{LyndonChristo}A word  $w$ is a Christoffel word if and only if  $w$ is a balanced Lyndon word. 
\elem

The next proposition follows from S\'e\'ebold, Richomme, Kassel and Reutenauer works \cite{ps1996,ps1998,gr2007,kr2007} and is proved in \cite{wfc1999}.

\bprop \label{christoMorph} Christoffel words and their conjugates are exactly the words obtained by the application of Sturmian morphisms to a letter. 
\eprop

{Lemma }\ref{LyndonChristo} {and Proposition} \ref{christoMorph} {have for consequence the following corollary.}

\bcor \label{agen} In the conjugation class of a Christoffel word, the Lyndon word is the Christoffel word. 
\ecor

{Note that Corollary }\ref{agen} { is the result we will extend as a definition of epichristoffel words.}


\begin{defn}\label{defEpiClass} A finite word  $w \in \A^*$ belongs to an  {\it  epichristoffel class} if it is the image of a letter by an episturmian morphism. 
\end{defn}

\begin{defn} A finite word $w \in \A^*$ is  {\it epichristoffel}\index{mot!epichristoffel@\'epichristoffel} if it is the unique Lyndon word occurring in an epichristoffel class.
\end{defn}

In the sequel, a word in an epichristoffel class will be called {\it $c$-epichristoffel}, {for short.}
{The following result insures that the epichristoffel classes are well-defined.} 

\bprop \label{propguil} Let $w$ and $w'$ be conjugate finite words. Then $w=\phi(u)$ and $w'=\phi'(u')$, with $\phi, \phi' \in \{\psi_a, \overline{\psi}_a \}$, for $u, u' \in \A^*$, $a \in \A$ if and only if $u$ and $u'$ are conjugate. 
\eprop

\Proof \begin{itemize} 
\item [($\Longrightarrow$)] Without loss of generality, we can suppose that  $\phi=\phi'=\psi_a$, since $\psi_a(w)=a\overline \psi_a(w) a^{-1}$ and so, $\psi_a(w)$ is conjugate to  $\overline \psi_a(w)$ for any word  $w$. Thus, we can write $w=a^{n_0}v[0]a^{n_1}v[1]\cdots a^{n_k}v[k]$, with $v[i] \neq a$ and $n_i >0$ for $0 \leq i \leq k$. Since $w=\psi_a(u)$, using injectivity of $\psi_a$, we have  $u=a^{n_0-1}v[0]a^{n_1-1}v[1]\cdots a^{n_k-1}v[k]$. Since $w$ and $w'$ are conjugate, we can write $w'=a^{\alpha}v[i]a^{n_{i+1}}v[i+1]\cdots a^{n_{i-1}}v[i-1]a^\beta$, with $\alpha+\beta = n_i$ and $\alpha \geq 1$. Thus, $u'=a^{\alpha-1}v[i]a^{n_{i+1}-1}v[i+1]\cdots a^{n_{i-1}-1}v[i-1]a^\beta$. Comparing $u$ and  $u'$, we conclude that  $u$ is conjugate to  $u'$. 
\item[($\Longleftarrow$)] If $u$ and $u'$ are conjugate, then there exist $v,t$ such that $u=vt$ and $u'=tv$. Applying respectively  the morphisms $\phi$ and $\phi'$ over $u$ and $u'$, we obtain $\phi(u)=\phi(v)\phi(t)$ and $\phi'(u')=\phi'(t)\phi'(v)$. If $\phi=\phi'$ the result follows. Otherwise, let us suppose $\phi=\psi_a$ and $\phi'=\overline \psi_a$. Then we conclude using the fact that $\psi_a(u)=a\overline \psi_a(u)a^{-1}$: 
$$\psi_a(u)=a\overline \psi_a(v)a^{-1}a\overline \psi(t)a^{-1}= a\overline \psi_a(v)\overline \psi_a(t)a^{-1}.$$
\vspace{-1.5cm}
 {\flushright \QED}
\end{itemize}

The finite factors of episturmian sequences, also called finite Arnoux-Rauzy words, have already been studied. In \cite {jp2002}, the authors used a subclass of $c$-epichristoffel words without mentioning that it is a generalization of Christoffel words. In their paper, they denoted by $h_n$, the standard episturmian words, that is the words obtained by the application of standard episturmian morphisms to a letter. The $c$-epichristoffel words are exactly  the set of all  conjugates of the standard episturmian words and the smallest one in the conjugacy class is epichristoffel. Notice that they form a subclass of the Arnoux-Rauzy word, since they all are factor of episturmian sequences, but any factor of episturmian sequence is not necessarily obtained by an episturmian morphism to a letter. For instance, the word $abacab\underline{aabac}ababacabaabacaba\cdots$ contains the finite Arnoux-Rauzy word  $aabac$ which is not $c$-epichristoffel.  

In \cite{jp2002}, the authors proved the $2$ following properties.

\bprop \textnormal{(\cite{jp2002}, prop. 2.8, prop. 2.12)} \label{propal1}Every standard episturmian word is primitive and can be written as the product of $2$ palindromic words.
\eprop

It is clear that any  standard episturmian word is conjugate to an epichristoffel word. Proposition \ref{propguil} can be used to show the converse.  Consequently, Proposition \ref{propal1} can be generalized for any $c$-epichristoffel word, {using the following lemma.}

\blem \textnormal{(\cite{djp2001}, Lemma 3)} \label{lemcon} {The word $u\in \A^*$ is a palindrome if and only if $\psi_a(u)a$ and $a \overline \psi_a(u)$ are so, $a \in \A$.}
\elem

\bprop  \label{propal} Every $c$-epichristoffel word is primitive and can be written as the product of $2$ palindromic words.
\eprop

\Proof {By induction over the number of morphisms. For a single morphism applied over a letter, we get  $w=ab$, with $a, b \in \A$ and $a\neq b$, which is the product of two palindromes. Let us suppose that for a $c$-epichristoffel word $w$, there exist palindromic words $u$, $v$ such that $w=uv$. Let $x=\psi _c(w)=\psi _c(uv)=\psi _c(u) \psi _c(v)$ (resp. $x=\overline \psi_c(w)=\overline \psi_c(u)\overline \psi_c(v)$), for $c \in \A$. Then $x=\psi _c(u) cc ^{-1} \psi _c(v)$ (resp. $x=\overline \psi _c(u) c ^{-1} c\overline \psi _c(v)$), where $\psi _c(u) c$, $c ^{-1} \psi _c(v)$ (resp. $\overline \psi _c(u) c ^{-1}$, $c \overline \psi _c(v)$) are palindromic words by Lemma} \ref{lemcon}.  \QED

Let now show how some of the properties of Christoffel words can be generalized to epichristoffel words.

Recall that for Christoffel {words}, we have:
\begin{theo} \label{thdldl}{\rm \cite{dldl2006}} Let $w$ be a non empty finite word. The following conditions are equivalent: \vspace{-0.2cm}
\begin{itemize}
\item [\rm i)] $w$ is a factor of a Sturmian sequence;
\item [\rm ii)] the fractionnary root $z_w$ of $w$ is conjugate to a Christoffel word.
\end{itemize}
\end{theo}

First, note that the equivalence in Theorem  \ref{thdldl} cannot be generalized to epichristoffel words. Indeed, let us consider the episturmian sequence
$$s=aabaacaabaacaabaabaa \cdot{ caabaacaabaaa}\cdots$$
Then  $w=caabaacaabaaa$ is a factor of  $s$, but its fractionnary root $z_w=w$ is not $c$-epichristoffel, as  we will see later in Example \ref{exa35}.

On the other hand, the converse holds for episturmian sequences and epichristoffel words.

\begin{theo} Let $w$ be a non empty word such that its fractionnary root is $c$-epichristoffel. Then $w$ is a factor of an episturmian sequence.
\end{theo}

\Proof Let $w=z^k_w$, with $k\geq 1 \in \Q$, $z_w$ the fractionnary root of $w$. Let us suppose that $z_w$ is $c$-epichristoffel. Thus there exist  $x \in \A^*$ and $a \in \A$ such that  $\phi^{(0)}\phi^{(1)}\cdots \phi^{(n)}(a)=z_w$, with $\phi^{(i)}\in \{\psi_{x[i]}, \overline \psi_{x[i]}\}$.  Then $w$ is a factor of  $z_w^{\lceil k \rceil}=(\phi^{(0)}\phi^{(1)}\cdots \phi^{(n)} (a))^{\lceil k \rceil}=\phi^{(0)}\phi^{(1)}\cdots \phi^{(n)}(a^{\lceil k \rceil})$. It is sufficient to take an episturmian sequence having $a^{\lceil k \rceil}$ as a factor  and apply the morphism $\phi^{(0)}\phi^{(1)}\cdots \phi^{(n)}$: we obtain that  $\phi^{(0)}\phi^{(1)}\cdots \phi^{(n)} (a^{\lceil k \rceil})$ is a factor of an episturmian sequence and so is $w$. \QED

\bprop Let $w \in \A^*$ be a $c$-epichristoffel word. Then, the set of factors of length $\leq |w|$ of its conjugacy class  is closed under mirror image.
\eprop

\Proof First note that the set of factors of length $\leq |w|$ of the epichristoffel class of $w$  is the same as the one of  $w^2$. Since any $c$-epichristoffel word  $w$ is the product of $2$ palindromes (by Proposition \ref{propal}), let $w=p_1p_2$, with $p_1$, $p_2$ palindromes. Then $w^2=p_1p_2p_1p_2$ and it follows that  $\widetilde{w}=\widetilde{p_1p_2}=p_2p_1$ is a factor of $w^2$. Thus, the mirror image of any factor of  $w$ is also a factor of  $w^2$ and consequently,  is in the epichristoffel class of  $w$. \QED

{\rem The right palindromic closure of a $c$-epichristoffel word is often {a} prefix of  $w^2$, but it is not the case in general. It suffices to take the word  $w=abcbab$ for which  $w^{(+)}=abcbab\cdot cba$. 
}

For Christoffel words, we have:
\blem {\rm \cite{dlm1994}}A Christoffel word can always be written as the product of two Christoffel words. 
\elem

But:
\blem \label{lemneg} An epichristoffel word cannot always be written as the product of two epichristoffel words. \elem

\Proof It is sufficient to consider the epichristoffel word $aabacab$. The only decompositions in $c$-epichristoffel factors are $a\cdot abacab$ and $aab\cdot acab$, but $abacab$ and $acab$ are not Lyndon words, assuming $a< b< c$. 
\QED

\blem Any $c$-epichristoffel word having length $> 1$ can be non-uniquely written as the product of two  $c$-epichristoffel words.
\elem

\Proof For the non unicity, it is sufficient to consider the example of the word  $aabacab$ given in the proof of  Lemma \ref{lemneg}. By definition, any $c$-epichristoffel word can be written as  $\phi^{(0)}\phi^{(1)} \cdots \phi^{(n-1)}(a)$, with $a \in \A$, $\phi^{(i)}\in \{\psi_{w[i]}, \overline \psi_{w[i]} \}$, $w \in \A^n$ and $w[n-1]\neq a$. Assume $\phi^{(n-1)}=\psi_{w[n-1]}$. To prove the existence of the product,  it is then  sufficient to consider the words $\phi^{(0)} \phi^{(1)}\cdots \phi^{(n-1)}(w[n-1])$ and $\phi^{(0)}\phi^{(1)}\cdots \phi^{(n-2)}(a)$, since
 \vspace{-0.2cm}
 \begin{eqnarray*}
\phi^{(0)}\phi^{(1)}\cdots \phi^{(n-1)}(a)&=& \phi^{(0)}\phi^{(1)}\cdots \phi^{(n-2)}(w[n-1]a)\\
&=& \phi^{(0)}\phi^{(1)}\cdots \phi^{(n-2)}(w[n-1]) \cdot \phi^{(0)}\phi^{(1)}\cdots \phi^{(n-2)}(a).
\end{eqnarray*} 
The case $\phi^{(n-1)}=\overline \psi_{w[n-1]}$ is analogue: we would have obtained a conjugate.
{\flushright \QED}
\vspace{-0.2cm}





\section{Epichristoffel $k$-tuples}

Recall from \cite{bl1993} that for a given $(p,q)$, with $p, q \in \N$, there exists a Christoffel word with occurrence numbers of letters  $p$ and $q$ if and only if  $p$ and $q$ are relatively primes. Moreover, it is possible to construct the corresponding Christoffel word, using a Cayley graph (see \cite{br2006}).

 In this section, we give an algorithm which determines if there exists or not an epichristoffel word  $w$ over the alphabet  $\A =\{a_0,a_1,\ldots,a_{k-1}\}$ such that  $p=(p_0,p_1,\dots, p_{k-1})$ with $p_i=|w|_{a_i}$, for $0\leq i \leq k-1$. If so, we also give an algorithm that constructs it.

\begin{defn} Let $p=(p_0,p_1,\ldots,p_{k-1})$ be a $k$-tuple of non negative integers. Then the  {\it operator} $T: \N^k \rightarrow \Z^k$ is defined over the $k$-tuple $p$ as
\vspace{-0.2cm}
$$T(p)=T(p_0,p_1,\ldots,p_{k-1})=(p_0,p_1,\ldots, p_{i-1},\left(p_i -  \sum _{j=0, j \neq i}^{k-1} p_j\right ), p_{i+1}, \ldots,p_{k-1}),$$ 
\vspace{-0.2cm}
where  $p_i \geq p_j$, $\forall j\neq i$.
\end{defn}

\bprop \label{ktuplets} Let $p$ be a $k$-tuple. There exists an epichristoffel word with occurrence numbers of letters $p$ if and only if  iterating $T$ over $p$ yields a $k$-tuple $p'$ with $p'_j=0$ for $j\neq m$ and $p'_m=1$, for a unique $m$ such that  $0\leq m \leq k-1$.
\eprop

The idea of using the operator $T$ comes from  the algorithm computing the greatest common divisor of $3$ integers as described in  \cite{cmr1999} and of the tuples described in \cite{jj2000}. 

{\exa \label{exa35} There is no epichristoffel word with the occurrence numbers of letters $(2,2,9)$. Indeed, $T(2,2,9)=(2,2,5)$, $T^2(2,2,9)=T(2,2,5)=(2,2,1)$, $T^3(2,2,9)=T(2,2,1)=(2,-1,1)$.
On the other hand, the $6$-tuple $q=(1,1,2,4,8,16)$ does so: 
\begin{eqnarray*}
T(1,1,2,4,8,16)&=&(1,1,2,4,8,0)\\
T^2(q)&=&T(1,1,2,4,8,0)=(1,1,2,4,0,0)\\
T^3(q)&=&T(1,1,2,4,0,0)=(1,1,2,0,0,0)\\
T^4(q)&=&T(1,1,2,0,0,0)=(1,1,0,0,0,0)\\
T^5(q)&=&T(1,1,0,0,0,0)=(1,0,0,0,0,0).
\end{eqnarray*}
}

Some lemmas are required in order to prove Proposition \ref{ktuplets}.

\blem \label{freq} Let $w=\phi(u)$, with $\phi \in \{\psi_{a_0}, \overline \psi_{a_0}\}$, $\A=\{a_0,a_1,\ldots, a_{k-1}\}$ and $u\in \A^*$. Then
\vspace{-0.2cm}
\begin{itemize}
\item [\rm i)] $\displaystyle |w|_{a_0}=\sum_{i=0}^{k-1} |u|_{a_i}=|u|$;
\vspace{-0.2cm}
\item [\rm ii)] $\displaystyle |w|_{a_0}=|u|_{a_0}+\sum_{i=1}^{k-1} |w|_{a_i}$.
\end{itemize}
\elem 

\Proof The first equality comes from the definition of $\psi_{a_0}$ and $\overline \psi_{a_0}$. For each letter $\alpha \neq a_0$, $\psi_{a_0}(\alpha)=a_0\alpha$,  $\overline \psi_{a_0}=\alpha a_0$ and $\overline\psi_{a_0}=\psi_{a_0}(a_0)=a_0$: $\phi$ adds as much $a_0$ as the occurrence numbers of the other letters in the word $u$. The second equality follows from the first one, since $|w|_{a_i}=|u|_{a_i}$ for $i\neq0$. 
\vspace{-0.8cm}
{\flushright \QED }

\blem \label{letterMax} Let $w\in \A^*$ be a $c$-epichristoffel word. Then, there exist a $c$-epichristoffel word  $u \in \A^*$, $|u| >1$ and an episturmian morphism $\phi \in \{\psi_{a_0}, \overline \psi_{a_0}\}$, with $a_0 \in \A$, such that $w=\phi(u)$ if and only if $|w|_{a_0} > |w|_{a_i}$ for all $a_i \in \A$, $i \neq 0$. 
\elem

\Proof 
\begin{itemize}
\item [($\Longrightarrow$)] By contradiction. Let us suppose there exists $u$ with $|u|>1$ such that $w=\phi(u)$ and $|w|_{a_0}$ is not maximum. Then, there exists at least one letter $a_i \in \A$ such that $|w|_{a_i} \geq |w|_{a_0}$. Without loss of generality, let us suppose that $i=1$. By Lemma \ref{freq},  $|w|_{a_0}=\sum_{i=0}^{k-1} |u|_{a_i}=|u|_{a_0}+|w|_{a_1}+\sum_{i =2}^{k-1}|u|_{a_i}$ that implies  $|w|_{a_0}-|w|_{a_1}=|u|_{a_0}+\sum_{i=2}^{k-1}|u|_{a_i} \leq 0$, which is possible only if  $|u|_{a_i}=0$ for all $i \neq 1$ and then $|w|_{a_1}=|w|_{a_0}$. Hence, we would have that  $u={a_1}^n$ and $w=\phi({a_1}^n)$.  The only possibility is that $n=1$, since a $c$-epichristoffel word is primitive. Then $|u|=1$: contradiction. Hence, if $w=\phi(u)$, with $|u|>1$,  $|w|_{a_0}$ is maximum.

\item [($\Longleftarrow$)] Let us now suppose that $|w|_{a_0} > |w|_{a_i}$ for all $a_i \in \A$, $i \neq 0$. Since $w$ is $c$-epichristoffel, there exist an episturmian morphism $\phi \in \{ \psi_{a_i}, \overline \psi_{a_i}\}$  and a $c$-epichristoffel word $u \in \A^*$ such that $\phi(u)=w$. Let us suppose that $i \neq 0$. Using Lemma \ref{freq}, $|w|_{a_i}=|w|_{a_0}+|u|_{a_i}+ \sum_{1\leq j \leq k-1, j\neq i}|w|_{a_j}$. Since $|w|_{a_0}> |w|_{a_i}$, it implies that $|u|_{a_i}+\sum_{1\leq j\leq k-1,j\neq i}|w|_{a_j}<0$, which is impossible. Thus $i=0$. \QED
%
\end{itemize}

An interesting consequence of Lemma \ref{letterMax} is the following.

\bprop \label{propUnique} Let $u$ and $v$ be $c$-epichristoffel words. If $|u|_\alpha=|v|_\alpha$ for all $\alpha \in \A$, then $u$ and $v$ are conjugate. In other words, a $k$-tuple of occurrence numbers of letters  determines at most one epichristoffel conjugacy class.
\eprop

\Proof {By induction. The result is true when $|u|=|v| \leq 2$. Assume by now that $|u| \geq 3$. By definition of epichristoffel words, there exist letters $a$ and $b$, and epichristoffel words $u', v'$ such that $u=\phi(u')$, $v=\phi'(v'), \phi \in \{\psi_a,\overline \psi_a\}$ and $\phi' \in \{\psi_b,\overline\psi_b\}$. From $|u| \geq 3$ and definitions  of morphisms $\psi_a, \overline \psi_a, \psi_b, \overline \psi_b$, we get $|u'| \geq 2, |v'| \geq 2$. From Lemma} \ref{letterMax} {and the fact that $|u|_\alpha=|v|_\alpha$ for all letters $\alpha$, it comes that $a=b$ (and $|u|_a=|v|_a \geq |u|_\alpha=|v|_\alpha$ for all letters $\alpha$). Now from definition of $u'$ and $v'$ and properties of $u$ and $v$, we deduce that $|u'|_\alpha=|v'|_\alpha$ for all letters $\alpha$. By inductive hypothesis, $u'$ and $v'$ are conjugate. Proposition} \ref{propguil} {allows to conclude.} \QED

%


The algorithm induced by the iteration of Lemma \ref{letterMax}  leads to a construction of words which are images of a letter by an episturmian morphism, that is $c$-epichristoffel words. Indeed, iterating $T$ gives a construction of an $c$-epichristoffel word with  $p$ describing the occurrence numbers of letters. We take $p$ as the initial $k$-tuple. The iteration over $p$ of the operator $T$ described previously yields a finite sequence of $k$-tuples $p^{(0)}$, $p^{(1)}$, $p^{(2)}, \dots$  We do as in Proposition  \ref{ktuplets}, applying  the operator $T$ and moreover, we keep an important information that allows us to construct the word: the letter with maximal number of occurrences.  Let 
\vspace{-0.2cm}
$$p^{(s)} \xrightarrow [ ]{\text{$i$}}  p^{(s+1)}$$
denote the relation $T(p^{(s)})=p^{(s+1)}$, where $p^{(s)}_i$ is the maximal integer of $p^{(s)}$. 

Then, performing  $T$ until  $p^{(r)}_i=0$ for all $i$ except for one  $i_{r-1}$ for which $p^{(r)}_{i_{r-1}}=1$, we get the  sequence of  $k$-tuples  
$$p^{(0)} \xrightarrow [ ]{\text{$i_0$}}  p^{(1)} \xrightarrow [ ]{\text{$i_1$}}  p^{(2)} \xrightarrow [ ]{\text{$i_2$}}  \cdots \xrightarrow [ ]{\text{$i_{r-2}$}}  p^{(r-1)} \xrightarrow [ ]{\text{$i_{r-1}$}}  p^{(r)}.$$ 
Then, 
$$\psi_{a_{i_0}}(\psi_{a_{i_1}}(\dots(\psi_{a_{i_{r-1}}}(\alpha))\dots))$$
is  a $c$-epichristoffel word having  $p$ as occurrence numbers of letters, with  $\alpha$ the letter such that $p^{(r)}_{i_{r-1}}=1$. The epichristoffel word is the Lyndon word of the conjugacy class of the word obtained. Here, Proposition \ref{propUnique} insures that it is sufficient to consider the standard episturmian morphism in order to construct a $c$-epichristoffel word with $p$ describing the occurrences of the letters.  \\

\noindent {\it {Proof of  Proposition} \ref{ktuplets}.}  Follows directly from Lemmas \ref{freq}, \ref{letterMax} and from the ideas described in the previous paragraph. The only difficulty concerns the last iteration, that is when $w=\phi(u)$, with $|w|_{a_0}$ not maximum. As seen in the previous proof, it implies that $u=a_1$ and $w=\phi(a_1) \in \{a_0a_1, a_1a_0\}$, which is clearly a $c$-epichristoffel word. Notice here that $\psi_{a_0}(a_1)=\overline \psi_{a_1}(a_0)$ and $\overline \psi_{a_0}(a_1)= \psi_{a_1}(a_0)$ are conjugate.    \QED

{\exa For the triplet $(5,10,16)$ describing the occurrence numbers of respectively the letters $a,b$ and $c$, the sequence obtained is
$$(5,10,16) \xrightarrow[ ]{\text{$c$}} (5,10,1) \xrightarrow[ ]{\text{$b$}} (5,4,1) \xrightarrow [ ]{\text{$a$}}(0,4,1)\xrightarrow[ ]{\text{$b$}} (0,3,1)\xrightarrow[ ]{\text{$b$}} (0,2,1)\xrightarrow[ ]{\text{$b$}} (0,1,1) \xrightarrow[]{\text{$b$}} (0,0,1).$$ 

Performing the algorithm, we find the word 
\vspace{-0.2cm}
\begin{eqnarray*}
\psi_{cbabbbb}(c)&=&\psi_{cbabbb}(bc)\\
&=&\psi_{cbabb}(\psi_b(bc))\\
&=&\psi_{cbab}(\psi_b(bbc))\\
&=&\psi_{cba}(\psi_b(bbbc))\\
&=&\psi_{cb}(\psi_a(bbbbc))\\
&=&\psi_{c}(\psi_b(ababababac))\\
&=&\psi_{c}(babbabbabbabbabc)\\
&=& cbcacbcbcacbcbcacbcbcacbcbcacbc.\\
\end{eqnarray*}
Since it is obtained by a standard episturmian morphism to a letter, this standard episturmian word is a representant of the epichristoffel conjugacy class. Moreover,  its conjugate which is a Lyndon word, and so, an epichristoffel word, is $acbcbcacbcbcacbcbcacbcbcacbc\cdot cbc$ for the order $a<b<c$.
}

Note that in the previous example, the choice of the last transition is arbitrary: we could have chosen the transition  $(0,1,1) \xrightarrow[]{\text{$c$}}(0,1,0)$ instead of $(0,1,1) \xrightarrow[]{\text{$b$}} (0,0,1)$ and we would have obtained a conjugate of  $\psi_{cbabbbb}(c)$ which is also $c$-epichristoffel.

\section{Criteria to be in an epichristoffel class} 

Let us recall a characterization of words in the conjugacy class of a Christoffel word. 

\begin{theo} \label{factSturm} {\rm \cite{dldl2006}} Let $w \in \A^*$ be a primitive word. Every conjugate $w' $ is a factor of a Sturmian sequence, not necessarily the same, if and only if  $w$ is conjugate to a Christoffel word.
\end{theo}

The goal of this section is to prove the following generalization of  Theorem \ref{factSturm}.

\begin{theo} \label{leth} Let $w$ be a finite primitive word different from a letter. Then there exists an episturmian sequence $z$ such that all the conjugates of $w$ are factors of  $z$  if and only if  $w$ is a $c$-epichristoffel word.
\end{theo}

Note that in order to generalize Theorem \ref{factSturm} to a $k$-letter alphabet, $k \geq 3$, an additional condition is necessary: the conjugates must be factor of the {\bf same} episturmian sequence. For example, every conjugates of the word $abc$ are factors of episturmian sequences, but  $abc$ is not a $c$-epichristoffel word, since $T(1,1,1)=(1,1,-1)$. 






Let us recall the following results of Justin and Pirillo that allow us to write any episturmian sequence as the image by an episturmian morphism of an other episturmian sequence.

\begin{cor} {\rm \cite{jp2002}}\label{cor37} Let $s \in \A^\omega$ and $\Delta=x[0]x[1]x[2]\cdots$, $x[i] \in \A$.  Then $s$ is a standard episturmian sequence with directive sequence $\Delta$ if and only if it exists an infinite sequence of sequences  $s^{(0)}=s, s^{(1)}, s^{(2)}, \ldots$ such that for any $i \in \N$, {$s^{(i-1)}=\psi_{x[i]}(s^{(i)})$}. 
\end{cor}

It can also be generalized to non standard episturmian sequences. In order to do so, let us recall what is  a  {\it spinned word}. Let $\overline \A= \{\overline a \, | \, a \in \A\}$. A letter $\overline x$ is considered as $x$ with {\it spin} $1$ while $x$ itself is considered as $x$ with spin $0$. Then, an {\it infinite spinned word} $\check s=\check s[0] \check s[1] \check s[2] \cdots $ is an element of $(\A \cup \overline \A)^\omega$. 

\begin{theo} {\rm \cite{jp2002}} \label{dirseq}A sequence $t\in \A^\omega$ is episturmian if and only if there exist a spinned sequence $\check{\Delta}=\check x[0]\check x[1]\check x[2]\cdots$, $\check x[i] \in \{\A \cup \overline \A\}$ and an infinite sequence of recurrent sequences $t^{(0)}=t$, $t^{(1)},t^{(2)}, \ldots$ such that for $i \in \N$, $t^{(i-1)}=\psi_{x[i]}(t^{(i)})$ if $\check x[i]$ has spin $0$ (resp. $\overline \psi_{x[i]}(t^{(i)})$ if $\check x[i]$ has spin $1$). Moreover $t$ is equivalent to the standard episturmian sequence with directive sequence $\Delta=x[0]x[1]\cdots$. 
\end{theo}

{Theorem} \ref{dirseq} {allows us to write the directive sequence of a non standard episturmian sequence, as we do in the following lemma.}

\blem \label{gp1} Let $\check \Delta(s)=(\check a)^k\check b \check z$ be the directive sequence of an episturmian sequence $s$, with $a \neq b \in \A$ and $z \in \A^\omega$. Then the blocks of $c\neq a$ have length $1$ and the blocks of $a$'s have  length  $\ell$, $k$ or $(k+1)$, where $\ell \leq k+1$ is the length of the block of $a$'s prefix of the sequence.
\elem

\Proof Let us consider the equivalent standard episturmian sequence $t$ directed by $\Delta(t)=a^kbz$. 
{By Corollary} \ref{cor37}, {$t=\psi_{a^kb}(t')$ for a standard episturmian word $t'$. Since $\psi_{a^kb}(a)=a^kba$, $\psi_{a^kb}(b)=a^kb$ and for $c \notin \{a,b\}$, $\psi_{a^kb}(c)=a^kba^kc$, the statement is true for $t$. Since the langage of $s$ and $t$ are equals, it only remains to consider the prefix of $s$ where a block of length $<k$ can appear. Indeed, for the episturmian sequence $s$, since it is directed by $\check \Delta (s)=(\check a)^k\check b \check z$, we easily deduce that $s$ begins by a prefix of $a$'s of length $\ell$ equals to the number of $\check a$ having spin $0$ in the prefix $(\check a)^k$ of its directive sequence, which is less or equal to $k$. \QED

}


{\rem An episturmian sequence may not have blocks of $a$'s of length $(k+1)$. It is the case if its directive sequence has the form $\check a^k\check z$, with $|\check z|_{\check a}=0$. 
}


One can be easily convinced of the following statement.

\blem \label{gp2} In an episturmian sequence $w=\psi_\alpha(t)$ or $w=\overline \psi _\alpha (t)$, any letter different from $\alpha$ is preceded and followed by the letter  $\alpha$, except for the first letter of the sequence, if it is different from $\alpha$. 
\elem



\blem \label{lemDeco} Let  $z=\psi _{a_0}(t)$ be a standard episturmian sequence and  $w=a_0ya_1$ a factor of $z$, with $a_0\neq a_1 \in \A$ and $y \in \A^*$. Then, there exists a factor $u$ of $t$ such that   $\psi_{a_0}(u)=w$. 
\elem

\Proof If $z=\psi_{a_0}(t)$, $t=t[0]t[1]t[2]\cdots$ and $\card(\A)=k$, then by the definition of $\psi$,  $z=\psi_{a_0}(t[0])\psi_{a_0}(t[1])\cdots \in \{a_0,a_0a_1, a_0a_2,\ldots, a_0a_{k-1}\}^\omega$. Since $w$ starts with $a_0$ and ends by $a_1$, then any  factor $w$ of $z=\psi_{a_0}(t)$ can be written as $w \in \{a_0,a_0a_1, a_0a_2, \ldots, a_0a_{k-1}\}^*$. Thus we can construct a word $u$ by associating to  $a_0a_i$ the letter  $a_i$ for $i\neq 0$ and to $a_0$ the letter  $a_0$. Thus,  $w$ is the image of the word  $u$ by the morphism $\psi_{a_0}$.  \QED

\bprop \label{prop4} Let $z=\psi _a(t)$, where $t$ and $z$ are standard episturmian sequences. Let $w$ be a factor of  $z$ not power of a letter, such that  $|w|>1$ and all its conjugates are also factors of  $z$. Then, there exists a factor $u$ of $t$ such that $w=\psi _a(u)$ or $w=\overline{\psi}_a(u)$. 
\eprop

\Proof Let $\beta, \gamma \in \A$, with $\beta,\gamma \neq a$, $y \in \A ^*$ and $w$ factor of $z$. There are $4$ cases to consider.
\begin{itemize}
\item [i)] $w=\beta y\gamma$:  its conjugate  $y\gamma \beta$ is not a factor of $z$, since any occurrence of the letter  $\beta$ is preceded by the letter $a$, by Lemma \ref{gp2}. Then $w$ does not satisfied the hypothesis.
\item [ii)] $w=ay\beta$: by Lemma \ref{lemDeco}, there exists  $u$ factor of  $t$ such that $\psi _a(u)=w$. 
\item [iii)] $w=\beta ya$: symmetric to the case  ii). If  $w=\beta ya$ is a factor of  $z=\psi _a(t)$ and satisfies the hypothesis, then there exists  $u$ factor of $t$ such that $\overline{\psi _a}(u)=w$.
\item [iv)] $w=aya$: rewrite  $w=a ^ my'a ^n$, with $m$, $n$ $\geq 1$ and $m, n$  maximum.  The factor $y'  $ is {not} empty, since $w$ is supposed not to be a power of a letter. Let us suppose that there exists $\beta \in \A$, $\beta \neq a$ such that  $w\beta=a ^my'a ^n\beta$ is a factor of  $z$. Since {by Lemma} \ref{gp1} any block of  $a$ has length  $k$ or $(k+1)$, {for some $k\in \N\setminus \{0\}$}, we have that  $n=k$ or $n=k+1$. On the other hand, by the hypothesis, the conjugate  $y'a ^{m+n} $  of $w=a ^m y' a ^n$ is also a factor of $z$. Thus  $m+n \leq k+1$. But since $m\neq 0$, the only possibility is that  $n=k$ and $m=1$. Consequently $w=ay'a ^k$. Its conjugate $y' a ^{k+1}$ is also a factor of  $z$ and since $y'$ does not start by  $a$ by the maximality of $m$, it should be preceded by $a$: $ay' a ^{k+1}=ay'a ^ka=wa$ is a factor of  $z$. Since $z$ is episturmian, $wa$ factor of $z$ implies that there exist $\ell \in \N$ and $\beta \neq a \in \A$ such that $wa^\ell\beta$ is so.  {By Lemma} \ref{lemDeco}, { there exists a word $u'=ua^{\ell-1}\beta$ such that $\psi_a(u')=wa^\ell\beta$. Since $\psi_a(a^{\ell-1}\beta)=a^\ell\beta$, $w=\psi_a(u)$. }
\end{itemize}
{\flushright \QED}



We can now prove our main Theorem.\\

\noindent {\it {Proof of Theorem} \ref{leth}}.  \begin{itemize}
\item [($\Longrightarrow$)] \begin{itemize} \item [i)] Let us suppose that all conjugates of $w$ are factor of a standard episturmian sequence $z=\psi_a(t)$. We proceed by induction on the number of morphisms. Since $z=\psi_a(t)$, by Proposition \ref{prop4}, there exists $u$ such that $w=\psi_a(u)$ or $w=\overline \psi_a(u)$. Let us now prove that all conjugates $u'$ of $u$ are also factors of $t$. Since $u, u'$ are conjugate, using Proposition \ref{propguil}, we have $\psi_a(u')$ is a conjugate of $ \psi_a(u)$. Hence, again by Proposition \ref{prop4}, there exists a factor $u''$ of $t$ with $\psi_a(u')=\psi_a(u'')$ or $ \psi_a(u')=\overline \psi_a(u'')$. The second case is possible only if $u''$ is a power of $a$ and then the first case holds. This first case  by injectivity of $\psi$ implies $u'=u''$, that is $u''$ is a factor of $t$. We then find a sequence of episturmian morphisms $\phi _0, \phi _1,.. ,\phi _k \in \{\psi_a,\overline \psi_a \, |\, a \in \A\}^{k+1}$ and a sequence of words  $w, w_1, w_2,... ,w_k$ such that $|w| \geq |w_1| \geq |w_2| \geq \ldots \geq |w_k|=1$, $w=\phi _0(\phi _1(...(\phi _k(w_k))...))$ and $w_i=\phi_i(\phi_{i+1}( \ldots (\phi_k(w_k))))$. Thus,  $w$ is the image of a letter by an episturmian morphism, implying that  $w$ is $c$-epichristoffel.
 \item [ii)] If $z$ is not standard, by  Definition \ref{episst}, we know that there exists an episturmian sequence $z'$ such that $F(z)=F(z')$. Thus, we can then consider  the sequence $z'$ and conclude as in i).
\end{itemize}

\item [($\Longleftarrow$)] Since $w$ is $c$-epichristoffel, we can write $w=f(a)$, where $f \in \mathscr E$ and $a \in \A$.  Let $s$ be an episturmian sequence having the factor $aa$ and let consider the episturmian sequence $f(s)$. Thus, it contains the factor $ww$ and we conclude.  

\end{itemize}
\vspace{-1cm}{\flushright \QED}

\section{Concluding remarks}
In this paper, we have most of the time consider the $c$-epichristoffel words, also known as the conjugates of the finite  standard episturmian words. Some of the properties of standard Sturmian words can be generalized naturally to the $c$-epichristoffel ones. We unfortunately didn't find a characterization of the epichristoffel word of each conjugacy class.  Geometrical properties of Christoffel words are well known and very interesting. It would be nice to know if there is a similar geometrical interpretation for the epichristoffel words.  In this paper, we only verify if a few properties of the Christoffel words could be generalized or not to the epichristoffel ones.  Since the literature of Christoffel words is wide, there are still a lot of open problems about epichristoffel words. For instance: do they satisfy a kind of balanced property? for a fixed $k \geq 3$, does there exist an epichristoffel word over a $k$-letter alphabet of any given length? is it possible to give a closed formula for the number of epichristoffel words of a given length? Episturmian morphisms have been extensively studied for instance in \cite{jj2001,jp2002,gr20032,gr2003,jp2004,jj2005,gr20072}. It might be useful to use their properties to work on the epichristoffel words.  

Epichristoffel words are still more interesting since they seem to be related to the Fraenkel conjecture. This conjecture states that for a finite $k$-letter alphabet, there exists a unique infinite word, up to letter permutation and conjugation, that is balanced and has pair-wise distinct letter frequencies. This unique word, if it exists, is conjectured to be periodic and can be written as $p^\omega$, with $p$ an epichristoffel word. Then, knowing more about epichristoffel words might help to prove the Fraenkel conjecture.

\section*{Acknowledgments} 
This paper is an extended version of a paper presented in Mons (Belgium) during the 12th Mons Theoretical Computer Science days \cite{gp2008}.  The author would also like to thank Christophe Reutenauer for giving her the idea of considering this interesting class of words, for useful discussions and remarks. Many thanks also to the two anonymous referees whose suggestions and constructive remarks helped to improve considerably the quality of the paper.

\bibliographystyle{alpha} 
{\footnotesize
\bibliography{biblio}
}

\end{document}